\documentclass[12pt,english]{article}
\usepackage{babel}
\usepackage[cp1251]{inputenc}
\usepackage{amsmath,amssymb}
\usepackage{epsfig}

\begin{document}

\newpage


\begin{center}
{\bf \Large
Estimates of solutions in a model

\vskip10pt

of antiviral immune response\footnote{The study was carried out
within the framework of the state contract of the Sobolev Institute of Mathematics 
(project no.~FWNF-2022-0008).}}
\end{center}

\begin{center}
{\bf M.A.~Skvortsova}

Sobolev Institute of Mathematics, Novosibirsk, Russia \\

sm-18-nsu@yandex.ru
\end{center}

{\bf Abstract.}
We consider a model of antiviral immune response proposed in the works of G.I.~Marchuk.
The model is described by a system of delay differential equations.
The asymptotic stability of a stationary solution to the system corresponding
to a completely healthy organism is studied.
Estimates of the attraction set of the given stationary solution are obtained
and estimates of solutions characterizing the stabilization rate at infinity are established.
The results are obtained using the Lyapunov–Krasovskii functional.
\\

\\




\section{Problem statement}

In the paper we consider one of the models of immunology
proposed in the works of G.I.~Marchuk~---
a model of antiviral immune response
(see~\cite{Marchuk1991}, subsection~3.2, p.~115--126).
The model is described by a system of delay  differential equations
and has the following form:
$$
\frac{d}{dt}x_1(t)=\nu x_9(t)+nb_{95}x_5(t) x_9(t)
-\gamma_{18}x_1(t) x_8(t)
$$
$$
-\gamma_{12} M x_1(t)
-\gamma_{19} C^{*} (1-x_9(t)-x_{10}(t))x_1(t),
$$
$$
\frac{d}{dt}x_2(t)=\gamma_{21} M x_1(t)-\alpha_2 x_2(t),
$$
$$
\frac{d}{dt}x_3(t)=b_{32} \left[ 
\rho_{32} \xi(x_{10}(t)) x_2(t-\tau_3) x_3(t-\tau_3)
-x_2(t) x_3(t) \right]
$$
$$
-b_3 x_2(t) x_3(t) x_5(t)+\alpha_3 (X_3^{*}-x_3(t)),
$$
$$
\frac{d}{dt}x_4(t)=b_{42} \left[ 
\rho_{42} \xi(x_{10}(t)) x_2(t-\tau_4) x_4(t-\tau_4)
-x_2(t) x_4(t) \right]
$$
$$
-b_4 x_2(t) x_4(t) x_6(t)+\alpha_4 (X_4^{*}-x_4(t)),
$$
$$
\frac{d}{dt}x_5(t)=b_5 \Bigl[ 
\rho_5 \xi(x_{10}(t)) x_2(t-\tau_5) x_3(t-\tau_5) x_5(t-\tau_5)
-x_2(t) x_3(t) x_5(t) \Bigr]
$$
$$
-b_{59} x_5(t) x_9(t)
+\alpha_5 (X_5^{*}-x_5(t)),
$$
$$
\frac{d}{dt}x_6(t)=b_6 \Bigl[ 
\rho_6 \xi(x_{10}(t)) x_2(t-\tau_6) x_4(t-\tau_6) x_6(t-\tau_6)
-x_2(t) x_4(t) x_6(t) \Bigr]
$$
$$
+\alpha_6 (X_6^{*}-x_6(t)),
$$
$$
\frac{d}{dt}x_7(t)=b_7 \rho_7 \xi(x_{10}(t))
x_2(t-\tau_7) x_4(t-\tau_7) x_6(t-\tau_7)
+\alpha_7 (X_7^{*}-x_7(t)),
$$
$$
\frac{d}{dt}x_8(t)=\rho_8 x_7(t)-\gamma_{81} x_1(t) x_8(t)
-\alpha_8 x_8(t),
$$
$$
\frac{d}{dt}x_9(t)=\sigma C^{*} x_1(t)(1-x_9(t)-x_{10}(t))
-b_{95}x_5(t) x_9(t)-b_{10} x_9(t),
$$
$$
\frac{d}{dt}x_{10}(t)=b_{95}x_5(t) x_9(t)+b_{10} x_9(t)-\alpha_{10} x_{10}(t).
\eqno (1.1)
$$
Here
$x_1(t)$
is the number of viruses freely circulating in the organism,
$x_2(t)$
is the number of stimulated macrophages,
$x_3(t)$
is the number of $T$-lymphocytes-helpers of cellular immunity,
$x_4(t)$
is the number of $T$-lymphocytes-helpers of humoral immunity,
$x_5(t)$
is the number of $T$-effector cells,
$x_6(t)$
is the number of $B$-lymphocytes,
$x_7(t)$
is the number of plasma cells,
$x_8(t)$
is the number of antibodies,
$x_9(t)$
is the number of target organ cells infected with viruses,
$x_{10}(t)$
is non-functioning part of the target organ affected by viruses
($x_{10}=0$
for an unaffected organ and
$x_{10}=1$
for a fully affected organ).
It is assumed that
$\xi(x_{10})$
is a non-increasing non-negative function defined on
$[0,1]$
and satisfying the Lipschitz condition on this segment,
which takes into account the disruption
of the normal functioning of the immune system
due to significant organ damage;
$0 \leq \xi(x_{10}) \leq 1$,
$\xi(0)=1$,
$\xi(1)=0$.
The values
$X_k^{*}>0$,
$k=3,\dots,7$,
denote the constant levels of the corresponding cells
in a completely healthy organism.
All coefficients of the system and delay parameters
are assumed to be constant and positive.
A more detailed description of the model is contained in~\cite{Marchuk1991}.

For system~(1.1), we set the initial conditions:
$$
x_j(t)=\varphi_j(t),
\quad
t \in [-\tau,0],
\quad
\tau=\max\limits_{k=3,\dots,7} \tau_k,
\quad
j=1,\dots,10,
\eqno (1.2)
$$
where
$\varphi_j(t)$,
$j=1,\dots,10$,
are given continuous functions, herewith
$$
\varphi_j(t) \ge 0,
\quad
t \in [-\tau,0],
\quad
j=1,\dots,10,
\quad
\max\limits_{t \in [-\tau,0]} \varphi_{10}(t)<1.
\eqno (1.3)
$$
It is well known that a solution to the initial value problem~(1.1)--(1.3)
exists, is unique, has non-negative components, and there exists a point
$t'>0$
such that
$x_{10}(t)<1$
for
$t\in (0,t')$.
Define the point
$t'$
as follows.
We put
$t'=\infty$,
if
$x_{10}(t)<1$
for
$t>0$,
and we put
$t'<\infty$
if
$x_{10}(t)<1$
for
$t \in (0,t')$
and
$x_{10}(t')=1$.
In the latter case, we assume that the solution to problem~(1.1)--(1.3)
is defined only for
$t \in (0,t')$.
From a biological point of view,
this means that over a period of time
$t \in (0,t')$,
viruses have completely infected the organ.
Below we will formulate the conditions
for the system parameters and for the initial data,
under which the inequality
$x_{10}(t)<1$
will hold for all
$t>0$.

As noted in~\cite{Marchuk1991}, system~(1.1) has a stationary solution
corresponding to a completely healthy organism:
$$
X^{*}=(X_1^{*},\dots,X_{10}^{*})^{\rm T}
=\left( 0,0,X_3^{*},X_4^{*},X_5^{*},X_6^{*},X_7^{*},
\frac{\rho_8}{\alpha_8} X_7^{*},0,0 \right)^{\rm T}.
\eqno (1.4)
$$
A sufficient condition for the asymptotic stability
of this stationary solution is the fulfillment of the inequality
(see~\cite{Marchuk1991}, p.~125)
$$
\Big(
\gamma_{12} M+\gamma_{18} \frac{\rho_8}{\alpha_8} X_7^{*}+\gamma_{19} C^{*}
\Big)
(b_{95} X_5^{*}+b_{10})
>\sigma C^{*} (\nu+nb_{95} X_5^{*}).
\eqno (1.5)
$$
The aim of this work is under condition~(1.5)
to indicate estimates of the attraction set of stationary solution~(1.4)
and to obtain estimates characterizing the rate of stabilization of solutions
to this stationary solution at infinity.
When obtaining the results, we will use
the method of Lyapunov~--Krasovskii functionals~\cite{DM2005}.



\section{The main result}

First, we reduce the problem of stability of the stationary solution~(1.4)
to the study of stability of the zero solution.
We change the variables
$$
x(t)=X^{*}+y(t),
$$
where
$x(t)=(x_1(t),\dots,x_{10}(t))^{\rm T}$,
$y(t)=(y_1(t),\dots,y_{10}(t))^{\rm T}$.
Then system~(1.1) is trans\-formed to the form
$$
\frac{d}{dt}y_1(t)=\nu y_9(t)+nb_{95} \Big( X_5^{*}+y_5(t) \Big) y_9(t)
-\gamma_{18}y_1(t) \Big( \frac{\rho_8}{\alpha_8} X_7^{*}+y_8(t) \Big)
$$
$$
-\gamma_{12} M y_1(t)
-\gamma_{19} C^{*} (1-y_9(t)-y_{10}(t))y_1(t),
$$
$$
\frac{d}{dt}y_2(t)=\gamma_{21} M y_1(t)-\alpha_2 y_2(t),
$$
$$
\frac{d}{dt}y_3(t)=b_{32} \left[ 
\rho_{32} \xi(y_{10}(t)) y_2(t-\tau_3) \Big( X_3^{*}+y_3(t-\tau_3) \Big)
-y_2(t) \Big( X_3^{*}+y_3(t) \Big) \right]
$$
$$
-b_3 y_2(t) \Big( X_3^{*}+y_3(t) \Big) \Big( X_5^{*}+y_5(t) \Big)-\alpha_3 y_3(t),
$$
$$
\frac{d}{dt}y_4(t)=b_{42} \left[ 
\rho_{42} \xi(y_{10}(t)) y_2(t-\tau_4) \Big( X_4^{*}+y_4(t-\tau_4) \Big)
-y_2(t) \Big( X_4^{*}+y_4(t) \Big) \right]
$$
$$
-b_4 y_2(t) \Big( X_4^{*}+y_4(t) \Big) \Big( X_6^{*}+y_6(t) \Big)-\alpha_4 y_4(t),
$$
$$
\frac{d}{dt}y_5(t)=b_5 \Bigl[ 
\rho_5 \xi(y_{10}(t)) y_2(t-\tau_5)
\Big( X_3^{*}+y_3(t-\tau_5) \Big)
\Big( X_5^{*}+y_5(t-\tau_5) \Big)
$$
$$
-y_2(t) \Big( X_3^{*}+y_3(t) \Big) \Big( X_5^{*}+y_5(t) \Big) \Bigr]
-b_{59} \Big( X_5^{*}+y_5(t) \Big) y_9(t)
-\alpha_5 y_5(t),
$$
$$
\frac{d}{dt}y_6(t)=b_6 \Bigl[ 
\rho_6 \xi(y_{10}(t)) y_2(t-\tau_6)
\Big( X_4^{*}+y_4(t-\tau_6) \Big)
\Big( X_6^{*}+y_6(t-\tau_6) \Big)
$$
$$
-y_2(t) \Big( X_4^{*}+y_4(t) \Big) \Big( X_6^{*}+y_6(t) \Big) \Bigr]
-\alpha_6 y_6(t),
$$
$$
\frac{d}{dt}y_7(t)=b_7 \rho_7 \xi(y_{10}(t))
y_2(t-\tau_7) \Big( X_4^{*}+y_4(t-\tau_7) \Big)
\Big( X_6^{*}+y_6(t-\tau_7) \Big)
-\alpha_7 y_7(t),
$$
$$
\frac{d}{dt}y_8(t)=\rho_8 y_7(t)
-\gamma_{81} y_1(t) \Big( \frac{\rho_8}{\alpha_8} X_7^{*}+y_8(t) \Big)
-\alpha_8 y_8(t),
$$
$$
\frac{d}{dt}y_9(t)=\sigma C^{*} y_1(t)(1-y_9(t)-y_{10}(t))
-b_{95} \Big( X_5^{*}+y_5(t) \Big) y_9(t)-b_{10} y_9(t),
$$
$$
\frac{d}{dt}y_{10}(t)=b_{95} \Big( X_5^{*}+y_5(t) \Big) y_9(t)
+b_{10} y_9(t)-\alpha_{10} y_{10}(t).
\eqno (2.1)
$$
The initial conditions have the form
$$
y_j(t)=\psi_j(t),
\quad
t \in [-\tau,0],
\quad
j=1,\dots,10,
\eqno (2.2)
$$
where
$\psi_j(t)$,
$j=1,\dots,10$,
are given continuous functions, herewith
$$
\psi_j(t) \ge -X_j^{*},
\quad
t \in [-\tau,0],
\quad
j=1,\dots,10,
\quad
\max\limits_{t \in [-\tau,0]} \psi_{10}(t)<1.
\eqno (2.3)
$$

{\bf Remark 1.}
Vector-functions
$\varphi(t)=(\varphi_1(t),\dots,\varphi_{10}(t))^{\rm T}$
and
$\psi(t)=(\psi_1(t),\dots, \psi_{10}(t))^{\rm T}$
are connected by the relation
$\psi(t)=\varphi(t)-X^{*}$.

Now we move on to the construction of the Lyapunov--Krasovskii functional. Let
$y(t)$
be a solution to the initial value problem~(2.1)--(2.3) defined for
$t \in (0,t')$.
We consider the functional
$$
V(t,y)=\sum\limits_{j=1}^{10} h_j y_j^2(t)
+\sum\limits_{k=3}^{7} \
\int\limits_{t-\tau_k}^{t} h_k \beta_k e^{-\ae_k (t-s)} y_2^2(s) ds.
\eqno (2.4)
$$
We define the values
$h_j>0$,
$j=1,\dots,10$,
$\beta_k>0$,
$\ae_k>0$,
$k=3,4,5,6,7$,
according to the following rule. Let
$\ae_k>0$,
$k=3,4,5,6,7$,
be arbitrary.
We also determine the values
$\theta_k>0$,
$k=3,4,5,6$,
in an arbitrary way.
We define the values
$\varepsilon$
and
$\delta$:
$$
\varepsilon=\frac{1}{2}
\left(
a_{11}+a_{99}-\sqrt{(a_{11}-a_{99})^2+4a_{19} a_{91}}
\right)>0,
\eqno (2.5)
$$
$$
0<\delta<\min\left\{
\varepsilon,
\alpha_2,
\alpha_3,
\alpha_4,
\alpha_5,
\alpha_6,
\alpha_7,
\alpha_8,
\alpha_{10}
\right\},
\eqno (2.6)
$$
where
$$
a_{11}=\Big(
\gamma_{12} M+\gamma_{18} \frac{\rho_8}{\alpha_8} X_7^{*}+\gamma_{19} C^{*}
\Big),
\quad
a_{99}=(b_{95} X_5^{*}+b_{10}),
\eqno (2.7)
$$
$$
a_{19}=(\nu+nb_{95} X_5^{*}),
\quad
a_{91}=\sigma C^{*}.
\eqno (2.8)
$$

{\bf Remark 2.}
The inequality
$\varepsilon>0$
follows from condition~(1.5) of the asymptotic stability
of the stationary solution~(1.4), which is equivalent to the inequality
$a_{11} a_{99}>a_{19} a_{91}$.

We set the values
$\beta_k>0$,
$k=3,4,5,6,7$:
$$
\beta_3=\frac{1}{\varepsilon_3}
b_{32} \rho_{32}
\Big( X_3^{*}+\theta_3 \Big)
e^{\ae_3 \tau_3 /2},
\eqno (2.9)
$$
$$
\beta_4=\frac{1}{\varepsilon_4}
b_{42} \rho_{42}
\Big( X_4^{*}+\theta_4 \Big)
e^{\ae_4 \tau_4 /2},
\eqno (2.10)
$$
$$
\beta_5=\frac{1}{\varepsilon_5}
b_5 \rho_5
\Big( X_3^{*}+\theta_3 \Big)
\Big( X_5^{*}+\theta_5 \Big)
e^{\ae_5 \tau_5 /2},
\eqno (2.11)
$$
$$
\beta_6=\frac{1}{\varepsilon_6}
b_6 \rho_6
\Big( X_4^{*}+\theta_4 \Big)
\Big( X_6^{*}+\theta_6 \Big)
e^{\ae_6 \tau_6 /2},
\eqno (2.12)
$$
$$
\beta_7=\frac{1}{\varepsilon_7}
b_7 \rho_7
\Big( X_4^{*}+\theta_4 \Big)
\Big( X_6^{*}+\theta_6 \Big)
e^{\ae_7 \tau_7 /2},
\eqno (2.13)
$$
where
$$
\varepsilon_3=\frac{2(\alpha_3-\delta)}
{\left(
b_{32} \rho_{32} \Big( X_3^{*}+\theta_3 \Big) e^{\ae_3 \tau_3 /2}
+\left( b_{32}+b_3 X_5^{*} \right) X_3^{*}
\right)},
\eqno (2.14)
$$
$$
\varepsilon_4=\frac{2(\alpha_4-\delta)}
{\left(
b_{42} \rho_{42} \Big( X_4^{*}+\theta_4 \Big) e^{\ae_4 \tau_4 /2}
+\left( b_{42}+b_4 X_6^{*} \right) X_4^{*}
\right)},
\eqno (2.15)
$$
$$
\varepsilon_5=\frac{(\alpha_5-\delta)}
{\left(
b_5 \rho_5
\Big( X_3^{*}+\theta_3 \Big)
\Big( X_5^{*}+\theta_5 \Big)
e^{\ae_5 \tau_5 /2}
+b_5 X_3^{*} X_5^{*}
\right)},
\eqno (2.16)
$$
$$
\varepsilon_6=\frac{2(\alpha_6-\delta)}
{\left(
b_6 \rho_6
\Big( X_4^{*}+\theta_4 \Big)
\Big( X_6^{*}+\theta_6 \Big)
e^{\ae_6 \tau_6 /2}
+b_6 X_4^{*} X_6^{*}
\right)},
\eqno (2.17)
$$
$$
\varepsilon_7=\frac{(\alpha_7-\delta)}
{b_7 \rho_7
\Big( X_4^{*}+\theta_4 \Big)
\Big( X_6^{*}+\theta_6 \Big)
e^{\ae_7 \tau_7 /2}}.
\eqno (2.18)
$$

We define the values
$h_j>0$,
$j=1,\dots,10$:
$$
h_1=\frac{1}{2(\varepsilon-\delta)}
\max\Bigg\{
\frac{a_{91}}{a_{19}}
\left[
\frac{(b_{59} X_5^{*})^2}
{\Big(
b_5 \rho_5
\Big( X_3^{*}+\theta_3 \Big)
\Big( X_5^{*}+\theta_5 \Big)
e^{\ae_5 \tau_5 /2}
+b_5 X_3^{*} X_5^{*}
\Big)^2}
+1
\right],
$$
$$
\left[
\frac{5(\gamma_{21} M)^2}{(\alpha_2-\delta)^2}
+\left( \gamma_{81} \frac{\rho_8}{\alpha_8} X_7^{*} \right)^2
\frac{(\alpha_7-\delta)^2}
{(b_7 \rho_7 \rho_8)^2
\Big( X_4^{*}+\theta_4 \Big)^2
\Big( X_6^{*}+\theta_6 \Big)^2
e^{\ae_7 \tau_7}}
\right]
\Bigg\},
\eqno (2.19)
$$
$$
h_2=\frac{5}{(\alpha_2-\delta)},
\eqno (2.20)
$$
$$
h_3=\frac{2(\alpha_3-\delta)}
{\left(
b_{32} \rho_{32}
\Big( X_3^{*}+\theta_3 \Big) e^{\ae_3 \tau_3 /2}
+\left( b_{32}+b_3 X_5^{*} \right) X_3^{*}
\right)^2},
\eqno (2.21)
$$
$$
h_4=\frac{2(\alpha_4-\delta)}
{\Big(
b_{42} \rho_{42}
\Big( X_4^{*}+\theta_4 \Big) e^{\ae_4 \tau_4 /2}
+\left( b_{42}+b_4 X_6^{*} \right) X_4^{*}
\Big)^2},
\eqno (2.22)
$$
$$
h_5=\frac{(\alpha_5-\delta)}
{\Big(
b_5 \rho_5
\Big( X_3^{*}+\theta_3 \Big)
\Big( X_5^{*}+\theta_5 \Big)
e^{\ae_5 \tau_5 /2}
+b_5 X_3^{*} X_5^{*}
\Big)^2},
\eqno (2.23)
$$
$$
h_6=\frac{2(\alpha_6-\delta)}
{\Big(
b_6 \rho_6
\Big( X_4^{*}+\theta_4 \Big)
\Big( X_6^{*}+\theta_6 \Big)
e^{\ae_6 \tau_6 /2}
+b_6 X_4^{*} X_6^{*}
\Big)^2},
\eqno (2.24)
$$
$$
h_7=\frac{(\alpha_7-\delta)}
{(b_7 \rho_7)^2
\Big( X_4^{*}+\theta_4 \Big)^2
\Big( X_6^{*}+\theta_6 \Big)^2
e^{\ae_7 \tau_7}},
\eqno (2.25)
$$
$$
h_8=\frac{(\alpha_7-\delta)^2 (\alpha_8-\delta)}
{(b_7 \rho_7 \rho_8)^2
\Big( X_4^{*}+\theta_4 \Big)^2
\Big( X_6^{*}+\theta_6 \Big)^2
e^{\ae_7 \tau_7}},
\eqno (2.26)
$$
$$
h_9=h_1 \frac{a_{19}}{a_{91}},
\eqno (2.27)
$$
$$
h_{10}=\frac{2(\alpha_{10}-\delta)}{(b_{95} X_5^{*}+b_{10})^2}.
\eqno (2.28)
$$

Functional
$V(t,y)$
is completely defined.
We also introduce the notation:
$$
\omega=\frac{1}{2} \min\left\{
2\delta,\ae_3,\ae_4,\ae_5,\ae_6,\ae_7 \right\},
\eqno (2.29)
$$
$$
q=2\Bigg(
\frac{nb_{95} \sqrt{h_1}}{\sqrt{h_5 h_9}}
+\frac{\gamma_{18}}{\sqrt{h_8}}
+\frac{\gamma_{19} C^{*}}{\sqrt{h_9}}
+\frac{\gamma_{19} C^{*}}{\sqrt{h_{10}}}
+\frac{b_3 X_3^{*} \sqrt{h_3}}{\sqrt{h_2 h_5}}
$$
$$
+\frac{b_4 X_4^{*} \sqrt{h_4}}{\sqrt{h_2 h_6}}
+\frac{b_5 X_5^{*} \sqrt{h_5}}{\sqrt{h_2 h_3}}
+\frac{b_6 X_6^{*} \sqrt{h_6}}{\sqrt{h_2 h_4}}
+\frac{b_{95}}{\sqrt{h_5}}
+\frac{b_{95} \sqrt{h_{10}}}{\sqrt{h_5 h_9}}
\Bigg).
\eqno (2.30)
$$

The following theorem is valid.

{\bf Theorem.}
{\it
Let condition~(1.5) be fulfilled.
Then for the components of the solution
to the initial value problem~(2.1)--(2.3),
with initial data satisfying the conditions
$$
\psi_j(t) \ge -X_j^{*},
\quad
t \in [-\tau,0],
\quad
j=1,\dots,10,
\quad
\sqrt{V(0,\psi)}<\frac{2\omega}{q},
\eqno (2.31)
$$
$$
\max\limits_{t \in [-\tau,0]} \psi_k(t) \leq \theta_k,
\quad
\frac{1}{\sqrt{h_k}} \frac{\sqrt{V(0,\psi)}}
{\displaystyle
\left( 1-\frac{q}{2\omega} \sqrt{V(0,\psi)} \right)} \leq \theta_k,
\quad k=3,4,5,6,
\eqno (2.32)
$$
$$
\max\limits_{t \in [-\tau,0]} \psi_{10}(t)<1,
\quad
\frac{1}{\sqrt{h_{10}}} \frac{\sqrt{V(0,\psi)}}
{\displaystyle
\left( 1-\frac{q}{2\omega} \sqrt{V(0,\psi)} \right)}<1,
\eqno (2.33)
$$
the estimates take place:
$$
y_j(t) \ge -X_j^{*},
\quad
j=1,\dots,10,
\quad
y_{10}(t)<1,
\quad t>0,
\eqno (2.34)
$$
$$
|y_j(t)| \leq \frac{1}{\sqrt{h_j}} \frac{\sqrt{V(0,\psi)}}
{\displaystyle
\left( 1-\frac{q}{2\omega} \sqrt{V(0,\psi)} \right)}
\, e^{-\omega t},
\quad
j=1,\dots,10,
\quad t>0.
\eqno (2.35)
$$
}

{\bf Remark 3.}
Inequalities~(2.31)--(2.33) are responsible
for the initial number of cells at which the organism recovers.
Estimates~(2.35) characterize the rate of recovery of the organism,
herewith the value
$e^{-\omega t}$
is responsible for the rate of recovery.
The inequality
$y_{10}(t)<1$
indicates that there is no complete organ lesion at any point in time.

{\bf Remark 4.}
Similar issues related to obtaining estimates of the attraction sets
of stationary solutions and estimates of solutions characterizing
the stabilization rate at infinity for some other models of immunology
were considered in~\cite{S2017}, \cite{S2018}.



\section{Proof of the theorem}



\subsection{Calculation of the derivative of the Lyapunov--Krasovskii
\linebreak
functional}

Let
$y(t)$
be a solution to the initial value problem~(2.1)--(2.3) defined for
$t \in (0,t')$.
We consider the Lyapunov--Krasovskii functional~(2.4):
$$
V(t,y)=\sum\limits_{j=1}^{10} h_j y_j^2(t)
+\sum\limits_{k=3}^{7} \
\int\limits_{t-\tau_k}^{t} h_k \beta_k e^{-\ae_k (t-s)} y_2^2(s) ds.
$$
We differentiate it along the solutions to system~(2.1):
$$
\frac{d}{dt} V(t,y)=\sum\limits_{j=1}^{10} 2h_j y_j(t) \frac{d}{dt} y_j(t)
+\sum\limits_{k=3}^{7} h_k \beta_k y_2^2(t)
-\sum\limits_{k=3}^{7} h_k \beta_k e^{-\ae_k \tau_k} y_2^2(t-\tau_k)
$$
$$
-\sum\limits_{k=3}^{7} \ae_k
\int\limits_{t-\tau_k}^{t} h_k \beta_k e^{-\ae_k (t-s)} y_2^2(s) ds.
$$
We take into account the fact that the inequalities are valid
$$
y_j(t) \ge 0, \quad j=1,2,9,10,
\quad
y_{10}(t) < 1,
\quad
\xi(y_{10}(t)) \leq 1,
\quad t \in (0,t').
$$
Then from system~(2.1) we get
$$
2h_1 y_1(t) \frac{d}{dt}y_1(t)
=-2h_1 \Big(
\gamma_{12} M+\gamma_{18} \frac{\rho_8}{\alpha_8} X_7^{*}+\gamma_{19} C^{*}
\Big) y_1^2(t)
$$
$$
+2h_1 (\nu+nb_{95} X_5^{*}) y_1(t) y_9(t)
+2h_1 nb_{95} y_1(t) y_5(t) y_9(t)
$$
$$
+2h_1 \Big(
-\gamma_{18} y_8(t)+\gamma_{19} C^{*} y_9(t)+\gamma_{19} C^{*} y_{10}(t)
\Big) y_1^2(t),
$$
\\
$$
2h_2 y_2(t) \frac{d}{dt}y_2(t)=-2h_2 \alpha_2 y_2^2(t)
+2h_2 \gamma_{21} M y_1(t) y_2(t),
$$
\\
$$
2h_3 y_3(t) \frac{d}{dt}y_3(t)
\leq
2h_3 b_{32} \rho_{32} y_2(t-\tau_3) \Big( X_3^{*}+y_3(t-\tau_3) \Big) y_3(t)
-2h_3 \alpha_3 y_3^2(t)
$$
$$
-2h_3 \left( b_{32}+b_3 X_5^{*} \right) X_3^{*} y_2(t) y_3(t)
-2h_3 b_3 X_3^{*} y_2(t) y_3(t) y_5(t),
$$
\\
$$
2h_4 y_4(t) \frac{d}{dt}y_4(t) \leq
2h_4 b_{42} \rho_{42} y_2(t-\tau_4) \Big( X_4^{*}+y_4(t-\tau_4) \Big) y_4(t)
-2h_4 \alpha_4 y_4^2(t)
$$
$$
-2h_4 \left( b_{42}+b_4 X_6^{*} \right) X_4^{*} y_2(t) y_4(t)
-2h_4 b_4 X_4^{*} y_2(t) y_4(t) y_6(t),
$$
\\
$$
2h_5 y_5(t) \frac{d}{dt}y_5(t) \leq
2h_5 b_5 \rho_5 y_2(t-\tau_5)
\Big( X_3^{*}+y_3(t-\tau_5) \Big)
\Big( X_5^{*}+y_5(t-\tau_5) \Big)
y_5(t)
$$
$$
-2h_5 \alpha_5 y_5^2(t)
-2h_5 b_5 X_3^{*} X_5^{*} y_2(t) y_5(t)
-2h_5 b_{59} X_5^{*} y_5(t) y_9(t)
-2h_5 b_5 X_5^{*} y_2(t) y_3(t) y_5(t),
$$
\\
$$
2h_6 y_6(t) \frac{d}{dt}y_6(t) \leq
2h_6 b_6 \rho_6 y_2(t-\tau_6)
\Big( X_4^{*}+y_4(t-\tau_6) \Big)
\Big( X_6^{*}+y_6(t-\tau_6) \Big)
y_6(t)
$$
$$
-2h_6 \alpha_6 y_6^2(t)
-2h_6 b_6 X_4^{*} X_6^{*} y_2(t) y_6(t)
-2h_6 b_6 X_6^{*} y_2(t) y_4(t) y_6(t),
$$
\\
$$
2h_7 y_7(t) \frac{d}{dt}y_7(t) \leq 2h_7 b_7 \rho_7
y_2(t-\tau_7) \Big( X_4^{*}+y_4(t-\tau_7) \Big)
\Big( X_6^{*}+y_6(t-\tau_7) \Big)
y_7(t)
$$
$$
-2h_7 \alpha_7 y_7^2(t),
$$
\\
$$
2h_8 y_8(t) \frac{d}{dt}y_8(t) \leq 2h_8 \rho_8 y_7(t) y_8(t)
-2h_8 \gamma_{81} \frac{\rho_8}{\alpha_8} X_7^{*} y_1(t) y_8(t)
-2h_8 \alpha_8 y_8^2(t),
$$
\\
$$
2h_9 y_9(t) \frac{d}{dt}y_9(t) \leq 2h_9 \sigma C^{*} y_1(t) y_9(t)
-2h_9 \left( b_{95} X_5^{*}+b_{10} \right) y_9^2(t)
-2h_9 b_{95} y_5(t) y_9^2(t),
$$
\\
$$
2h_{10} y_{10}(t) \frac{d}{dt}y_{10}(t)
=2h_{10} (b_{95} X_5^{*}+b_{10}) y_9(t) y_{10}(t)
+2h_{10} b_{95} y_5(t) y_9(t) y_{10}(t)
$$
$$
-2h_{10} \alpha_{10} y_{10}^2(t).
$$

Consequently, the inequality takes place
$$
\frac{d}{dt} V(t,y) \leq U_{01}+U_{02}+U_{03}+U_{\tau}
-\sum\limits_{k=3}^{7} \ae_k
\int\limits_{t-\tau_k}^{t} h_k \beta_k e^{-\ae_k (t-s)} y_2^2(s) ds,
\eqno (3.1)
$$
where
$$
U_{01}=-2h_1 \Big(
\gamma_{12} M+\gamma_{18} \frac{\rho_8}{\alpha_8} X_7^{*}+\gamma_{19} C^{*}
\Big) y_1^2(t)
-\left( 2h_2 \alpha_2-\sum\limits_{k=3}^{7} h_k \beta_k \right) y_2^2(t)
$$
$$
-2h_3 \alpha_3 y_3^2(t)
-2h_4 \alpha_4 y_4^2(t)
-2h_5 \alpha_5 y_5^2(t)
-2h_6 \alpha_6 y_6^2(t)
-2h_7 \alpha_7 y_7^2(t)
$$
$$
-2h_8 \alpha_8 y_8^2(t)
-2h_9 \left( b_{95} X_5^{*}+b_{10} \right) y_9^2(t)
-2h_{10} \alpha_{10} y_{10}^2(t),
\eqno (3.2)
$$
\\
$$
U_{02}=2h_1 (\nu+nb_{95} X_5^{*}) y_1(t) y_9(t)
+2h_2 \gamma_{21} M y_1(t) y_2(t)
$$
$$
-2h_3 \left( b_{32}+b_3 X_5^{*} \right) X_3^{*} y_2(t) y_3(t)
-2h_4 \left( b_{42}+b_4 X_6^{*} \right) X_4^{*} y_2(t) y_4(t)
$$
$$
-2h_5 b_5 X_3^{*} X_5^{*} y_2(t) y_5(t)
-2h_5 b_{59} X_5^{*} y_5(t) y_9(t)
-2h_6 b_6 X_4^{*} X_6^{*} y_2(t) y_6(t)
$$
$$
+2h_8 \rho_8 y_7(t) y_8(t)
-2h_8 \gamma_{81} \frac{\rho_8}{\alpha_8} X_7^{*} y_1(t) y_8(t)
+2h_9 \sigma C^{*} y_1(t) y_9(t)
$$
$$
+2h_{10} (b_{95} X_5^{*}+b_{10}) y_9(t) y_{10}(t),
\eqno (3.3)
$$
\\
$$
U_{03}=2h_1 nb_{95} y_1(t) y_5(t) y_9(t)
+2h_1 \Big(
-\gamma_{18} y_8(t)+\gamma_{19} C^{*} y_9(t)+\gamma_{19} C^{*} y_{10}(t)
\Big) y_1^2(t)
$$
$$
-2h_3 b_3 X_3^{*} y_2(t) y_3(t) y_5(t)
-2h_4 b_4 X_4^{*} y_2(t) y_4(t) y_6(t)
-2h_5 b_5 X_5^{*} y_2(t) y_3(t) y_5(t)
$$
$$
-2h_6 b_6 X_6^{*} y_2(t) y_4(t) y_6(t)
-2h_9 b_{95} y_5(t) y_9^2(t)
+2h_{10} b_{95} y_5(t) y_9(t) y_{10}(t),
\eqno (3.4)
$$
\\
$$
U_{\tau}
=2h_3 b_{32} \rho_{32} y_2(t-\tau_3) \Big( X_3^{*}+y_3(t-\tau_3) \Big) y_3(t)
$$
$$
+2h_4 b_{42} \rho_{42} y_2(t-\tau_4) \Big( X_4^{*}+y_4(t-\tau_4) \Big) y_4(t)
$$
$$
+2h_5 b_5 \rho_5 y_2(t-\tau_5)
\Big( X_3^{*}+y_3(t-\tau_5) \Big)
\Big( X_5^{*}+y_5(t-\tau_5) \Big)
y_5(t)
$$
$$
+2h_6 b_6 \rho_6 y_2(t-\tau_6)
\Big( X_4^{*}+y_4(t-\tau_6) \Big)
\Big( X_6^{*}+y_6(t-\tau_6) \Big)
y_6(t)
$$
$$
+2h_7 b_7 \rho_7
y_2(t-\tau_7) \Big( X_4^{*}+y_4(t-\tau_7) \Big)
\Big( X_6^{*}+y_6(t-\tau_7) \Big)
y_7(t)
$$
$$
-\sum\limits_{k=3}^{7} h_k \beta_k e^{-\ae_k \tau_k} y_2^2(t-\tau_k).
\eqno (3.5)
$$



\subsection{An estimate for $U_{03}$}

Consider function
$U_{03}$
defined in~(3.4).
We have
$$
U_{03} \leq 2h_1 nb_{95} \frac{V^{3/2}(t,y)}{\sqrt{h_1 h_5 h_9}}
+2 \left(
\gamma_{18} \frac{V^{1/2}(t,y)}{\sqrt{h_8}}
+\gamma_{19} C^{*} \frac{V^{1/2}(t,y)}{\sqrt{h_9}}
+\gamma_{19} C^{*} \frac{V^{1/2}(t,y)}{\sqrt{h_{10}}}
\right) V(t,y)
$$
$$
+2h_3 b_3 X_3^{*} \frac{V^{3/2}(t,y)}{\sqrt{h_2 h_3 h_5}}
+2h_4 b_4 X_4^{*} \frac{V^{3/2}(t,y)}{\sqrt{h_2 h_4 h_6}}
+2h_5 b_5 X_5^{*} \frac{V^{3/2}(t,y)}{\sqrt{h_2 h_3 h_5}}
$$
$$
+2h_6 b_6 X_6^{*} \frac{V^{3/2}(t,y)}{\sqrt{h_2 h_4 h_6}}
+2 b_{95} \frac{V^{3/2}(t,y)}{\sqrt{h_5}}
+2h_{10} b_{95} \frac{V^{3/2}(t,y)}{\sqrt{h_5 h_9 h_{10}}}.
$$
Therefore,
$$
U_{03} \leq q V^{3/2}(t,y),
\eqno (3.6)
$$
where
$q$
is defined in~(2.30).



\subsection{An estimate for $U_{\tau}$}

Consider function
$U_{\tau}$
defined in~(3.5).
We have
$$
U_{\tau}=2h_3 b_{32} \rho_{32}
y_2(t-\tau_3) \Big( X_3^{*}+y_3(t-\tau_3) \Big) y_3(t)
-h_3 \beta_3 e^{-\ae_3 \tau_3} y_2^2(t-\tau_3)
$$
$$
+2h_4 b_{42} \rho_{42}
y_2(t-\tau_4) \Big( X_4^{*}+y_4(t-\tau_4) \Big) y_4(t)
-h_4 \beta_4 e^{-\ae_4 \tau_4} y_2^2(t-\tau_4)
$$
$$
+2h_5 b_5 \rho_5
y_2(t-\tau_5) \Big( X_3^{*}+y_3(t-\tau_5) \Big)
\Big( X_5^{*}+y_5(t-\tau_5) \Big) y_5(t)
-h_5 \beta_5 e^{-\ae_5 \tau_5} y_2^2(t-\tau_5)
$$
$$
+2h_6 b_6 \rho_6
y_2(t-\tau_6) \Big( X_4^{*}+y_4(t-\tau_6) \Big)
\Big( X_6^{*}+y_6(t-\tau_6) \Big) y_6(t)
-h_6 \beta_6 e^{-\ae_6 \tau_6} y_2^2(t-\tau_6)
$$
$$
+2h_7 b_7 \rho_7
y_2(t-\tau_7) \Big( X_4^{*}+y_4(t-\tau_7) \Big)
\Big( X_6^{*}+y_6(t-\tau_7) \Big) y_7(t)
-h_7 \beta_7 e^{-\ae_7 \tau_7} y_2^2(t-\tau_7).
$$
Therefore,
$$
U_{\tau} \leq W_{\tau},
$$
where
$$
W_{\tau}=h_3 \frac{e^{\ae_3 \tau_3}}{\beta_3}
(b_{32} \rho_{32})^2
\Big( X_3^{*}+y_3(t-\tau_3) \Big)^2 y_3^2(t)
$$
$$
+h_4 \frac{e^{\ae_4 \tau_4}}{\beta_4}
(b_{42} \rho_{42})^2
\Big( X_4^{*}+y_4(t-\tau_4) \Big)^2 y_4^2(t)
$$
$$
+h_5 \frac{e^{\ae_5 \tau_5}}{\beta_5}
(b_5 \rho_5)^2
\Big( X_3^{*}+y_3(t-\tau_5) \Big)^2
\Big( X_5^{*}+y_5(t-\tau_5) \Big)^2 y_5^2(t)
$$
$$
+h_6 \frac{e^{\ae_6 \tau_6}}{\beta_6}
(b_6 \rho_6)^2
\Big( X_4^{*}+y_4(t-\tau_6) \Big)^2
\Big( X_6^{*}+y_6(t-\tau_6) \Big)^2 y_6^2(t)
$$
$$
+h_7 \frac{e^{\ae_7 \tau_7}}{\beta_7}
(b_7 \rho_7)^2
\Big( X_4^{*}+y_4(t-\tau_7) \Big)^2
\Big( X_6^{*}+y_6(t-\tau_7) \Big)^2 y_7^2(t).
$$
We represent
$W_{\tau}$
in the form
$$
W_{\tau}=W_{01}+R_{\tau},
$$
where
$$
W_{01}=h_3 \frac{e^{\ae_3 \tau_3}}{\beta_3}
(b_{32} \rho_{32})^2
\Big( X_3^{*}+\theta_3 \Big)^2 y_3^2(t)
$$
$$
+h_4 \frac{e^{\ae_4 \tau_4}}{\beta_4}
(b_{42} \rho_{42})^2
\Big( X_4^{*}+\theta_4 \Big)^2 y_4^2(t)
$$
$$
+h_5 \frac{e^{\ae_5 \tau_5}}{\beta_5}
(b_5 \rho_5)^2
\Big( X_3^{*}+\theta_3 \Big)^2
\Big( X_5^{*}+\theta_5 \Big)^2 y_5^2(t)
$$
$$
+h_6 \frac{e^{\ae_6 \tau_6}}{\beta_6}
(b_6 \rho_6)^2
\Big( X_4^{*}+\theta_4 \Big)^2
\Big( X_6^{*}+\theta_6 \Big)^2 y_6^2(t)
$$
$$
+h_7 \frac{e^{\ae_7 \tau_7}}{\beta_7}
(b_7 \rho_7)^2
\Big( X_4^{*}+\theta_4 \Big)^2
\Big( X_6^{*}+\theta_6 \Big)^2 y_7^2(t),
\eqno (3.7)
$$
\\
$$
R_{\tau}=h_3 \frac{e^{\ae_3 \tau_3}}{\beta_3}
(b_{32} \rho_{32})^2
\Big[
\Big( X_3^{*}+y_3(t-\tau_3) \Big)^2
-\Big( X_3^{*}+\theta_3 \Big)^2
\Big]
y_3^2(t)
$$
$$
+h_4 \frac{e^{\ae_4 \tau_4}}{\beta_4}
(b_{42} \rho_{42})^2
\Big[
\Big( X_4^{*}+y_4(t-\tau_4) \Big)^2
-\Big( X_4^{*}+\theta_4 \Big)^2
\Big]
y_4^2(t)
$$
$$
+h_5 \frac{e^{\ae_5 \tau_5}}{\beta_5}
(b_5 \rho_5)^2
\Big[
\Big( X_3^{*}+y_3(t-\tau_5) \Big)^2
\Big( X_5^{*}+y_5(t-\tau_5) \Big)^2
-\Big( X_3^{*}+\theta_3 \Big)^2
\Big( X_5^{*}+\theta_5 \Big)^2
\Big]
y_5^2(t)
$$
$$
+h_6 \frac{e^{\ae_6 \tau_6}}{\beta_6}
(b_6 \rho_6)^2
\Big[
\Big( X_4^{*}+y_4(t-\tau_6) \Big)^2
\Big( X_6^{*}+y_6(t-\tau_6) \Big)^2
-\Big( X_4^{*}+\theta_4 \Big)^2
\Big( X_6^{*}+\theta_6 \Big)^2
\Big]
y_6^2(t)
$$
$$
+h_7 \frac{e^{\ae_7 \tau_7}}{\beta_7}
(b_7 \rho_7)^2
\Big[
\Big( X_4^{*}+y_4(t-\tau_7) \Big)^2
\Big( X_6^{*}+y_6(t-\tau_7) \Big)^2
-\Big( X_4^{*}+\theta_4 \Big)^2
\Big( X_6^{*}+\theta_6 \Big)^2
\Big]
y_7^2(t),
\eqno (3.8)
$$
$\theta_k>0$,
$k=3,4,5,6$,
are arbitrary positive numbers. So,
$$
U_{\tau} \leq W_{01}+R_{\tau}.
\eqno (3.9)
$$



\subsection{An estimate for $U_{02}$}

Consider function
$U_{02}$
defined in~(3.3).
We have
$$
U_{02}=2h_1 (\nu+nb_{95} X_5^{*}) y_1(t) y_9(t)
+2h_2 \gamma_{21} M y_1(t) y_2(t)
$$
$$
-2h_3 \left( b_{32}+b_3 X_5^{*} \right) X_3^{*} y_2(t) y_3(t)
-2h_4 \left( b_{42}+b_4 X_6^{*} \right) X_4^{*} y_2(t) y_4(t)
$$
$$
-2h_5 b_5 X_3^{*} X_5^{*} y_2(t) y_5(t)
-2h_5 b_{59} X_5^{*} y_5(t) y_9(t)
-2h_6 b_6 X_4^{*} X_6^{*} y_2(t) y_6(t)
$$
$$
+2h_8 \rho_8 y_7(t) y_8(t)
-2h_8 \gamma_{81} \frac{\rho_8}{\alpha_8} X_7^{*} y_1(t) y_8(t)
+2h_9 \sigma C^{*} y_1(t) y_9(t)
$$
$$
+2h_{10} (b_{95} X_5^{*}+b_{10}) y_9(t) y_{10}(t).
$$
Hence the inequality follows
$$
U_{02} \leq W_{02},
\eqno (3.10)
$$
where
$$
W_{02}=h_1 (\nu+nb_{95} X_5^{*})
\left( \varepsilon_1 y_1^2(t)+\frac{y_9^2(t)}{\varepsilon_1} \right)
+h_2 \gamma_{21} M
\left( \varepsilon_2 y_2^2(t)+\frac{y_1^2(t)}{\varepsilon_2} \right)
$$
$$
+h_3 \left( b_{32}+b_3 X_5^{*} \right) X_3^{*}
\left( \varepsilon_3 y_3^2(t)+\frac{y_2^2(t)}{\varepsilon_3} \right)
+h_4 \left( b_{42}+b_4 X_6^{*} \right) X_4^{*}
\left( \varepsilon_4 y_4^2(t)+\frac{y_2^2(t)}{\varepsilon_4} \right)
$$
$$
+h_5 b_5 X_3^{*} X_5^{*}
\left( \varepsilon_5 y_5^2(t)+\frac{y_2^2(t)}{\varepsilon_5} \right)
+h_5 b_{59} X_5^{*}
\left( \varepsilon_{59} y_5^2(t)+\frac{y_9^2(t)}{\varepsilon_{59}} \right)
$$
$$
+h_6 b_6 X_4^{*} X_6^{*}
\left( \varepsilon_6 y_6^2(t)+\frac{y_2^2(t)}{\varepsilon_6} \right)
+h_8 \rho_8
\left( \varepsilon_{87} y_8^2(t)+\frac{y_7^2(t)}{\varepsilon_{87}} \right)
$$
$$
+h_8 \gamma_{81} \frac{\rho_8}{\alpha_8} X_7^{*}
\left( \varepsilon_{81} y_8^2(t)+\frac{y_1^2(t)}{\varepsilon_{81}} \right)
+h_9 \sigma C^{*}
\left( \varepsilon_9 y_9^2(t)+\frac{y_1^2(t)}{\varepsilon_9} \right)
$$
$$
+h_{10} (b_{95} X_5^{*}+b_{10})
\left( \varepsilon_{10} y_{10}^2(t)+\frac{y_9^2(t)}{\varepsilon_{10}} \right),
\eqno (3.11)
$$
where
$\varepsilon_3$,
$\varepsilon_4$,
$\varepsilon_5$,
$\varepsilon_6>0$
are defined in~(2.14)--(2.17),
$\varepsilon_1$,
$\varepsilon_2$,
$\varepsilon_{59}$,
$\varepsilon_{87}$,
$\varepsilon_{81}$,
$\varepsilon_9$,
$\varepsilon_{10}>0$
will be defined below.



\subsection{An estimate for $U_{01}+U_{02}+U_{\tau}$}

We estimate
$U_{01}+U_{02}+U_{\tau}$.
From inequalities~(3.9) and~(3.10), the estimate follows
$$
U_{01}+U_{02}+U_{\tau} \leq U_{01}+W_{01}+W_{02}+R_{\tau},
\eqno (3.12)
$$
where
$U_{01}$
is defined in~(3.2),
$W_{01}$
is defined in~(3.7),
$W_{02}$
is defined in~(3.11),
$R_{\tau}$
is defined in~(3.8).
We have
$$
U_{01}+W_{01}+W_{02}
=-\sum\limits_{j=1}^{10} h_j p_j y_j^2(t),
\eqno (3.13)
$$
where
$$
p_1=2 \Big(
\gamma_{12} M+\gamma_{18} \frac{\rho_8}{\alpha_8} X_7^{*}+\gamma_{19} C^{*}
\Big)
-\varepsilon_1 (\nu+nb_{95} X_5^{*})
$$
$$
-\frac{1}{\varepsilon_9} \sigma C^{*} \left( \frac{h_9}{h_1} \right)
-\frac{1}{\varepsilon_2} \gamma_{21} M \left( \frac{h_2}{h_1} \right)
-\frac{1}{\varepsilon_{81}}
\gamma_{81} \frac{\rho_8}{\alpha_8} X_7^{*}
\left( \frac{h_8}{h_1} \right),
$$
$$
p_2=2\alpha_2-\varepsilon_2 \gamma_{21} M
-\sum\limits_{k=3}^{7} \beta_k \left( \frac{h_k}{h_2} \right)
-\frac{1}{\varepsilon_3}
\left( b_{32}+b_3 X_5^{*} \right) X_3^{*}
\left( \frac{h_3}{h_2} \right)
$$
$$
-\frac{1}{\varepsilon_4} \left( b_{42}+b_4 X_6^{*} \right) X_4^{*}
\left( \frac{h_4}{h_2} \right)
-\frac{1}{\varepsilon_5} b_5 X_3^{*} X_5^{*}
\left( \frac{h_5}{h_2} \right)
-\frac{1}{\varepsilon_6} b_6 X_4^{*} X_6^{*}
\left( \frac{h_6}{h_2} \right),
$$
$$
p_3=2\alpha_3
-\frac{e^{\ae_3 \tau_3}}{\beta_3}
(b_{32} \rho_{32})^2 \Big( X_3^{*}+\theta_3 \Big)^2
-\varepsilon_3 \left( b_{32}+b_3 X_5^{*} \right) X_3^{*},
$$
$$
p_4=2\alpha_4
-\frac{e^{\ae_4 \tau_4}}{\beta_4}
(b_{42} \rho_{42})^2 \Big( X_4^{*}+\theta_4 \Big)^2
-\varepsilon_4 \left( b_{42}+b_4 X_6^{*} \right) X_4^{*},
$$
$$
p_5=2\alpha_5
-\frac{e^{\ae_5 \tau_5}}{\beta_5}
(b_5 \rho_5)^2
\Big( X_3^{*}+\theta_3 \Big)^2
\Big( X_5^{*}+\theta_5 \Big)^2
-\varepsilon_5 b_5 X_3^{*} X_5^{*}
-\varepsilon_{59} b_{59} X_5^{*},
$$
$$
p_6=2\alpha_6
-\frac{e^{\ae_6 \tau_6}}{\beta_6}
(b_6 \rho_6)^2
\Big( X_4^{*}+\theta_4 \Big)^2
\Big( X_6^{*}+\theta_6 \Big)^2
-\varepsilon_6 b_6 X_4^{*} X_6^{*},
$$
$$
p_7=2\alpha_7
-\frac{e^{\ae_7 \tau_7}}{\beta_7}
(b_7 \rho_7)^2
\Big( X_4^{*}+\theta_4 \Big)^2
\Big( X_6^{*}+\theta_6 \Big)^2
-\frac{1}{\varepsilon_{87}} \rho_8
\left( \frac{h_8}{h_7} \right),
$$
$$
p_8=2\alpha_8
-\varepsilon_{87} \rho_8
-\varepsilon_{81} \gamma_{81} \frac{\rho_8}{\alpha_8} X_7^{*},
$$
$$
p_9=2 \left( b_{95} X_5^{*}+b_{10} \right)
-\varepsilon_9 \sigma C^{*}
-\frac{1}{\varepsilon_1} (\nu+nb_{95} X_5^{*})
\left( \frac{h_1}{h_9} \right)
$$
$$
-\frac{1}{\varepsilon_{59}} b_{59} X_5^{*}
\left( \frac{h_5}{h_9} \right)
-\frac{1}{\varepsilon_{10}} (b_{95} X_5^{*}+b_{10})
\left( \frac{h_{10}}{h_9} \right),
$$
$$
p_{10}=2\alpha_{10}
-\varepsilon_{10} (b_{95} X_5^{*}+b_{10}).
$$

Taking into account the notation~(2.7), (2.8) and definition~(2.27),
we have
$$
p_1=2 a_{11}-a_{19} \left( \varepsilon_1+\frac{1}{\varepsilon_9} \right)
-\frac{1}{\varepsilon_2} \gamma_{21} M \left( \frac{h_2}{h_1} \right)
-\frac{1}{\varepsilon_{81}}
\gamma_{81} \frac{\rho_8}{\alpha_8} X_7^{*}
\left( \frac{h_8}{h_1} \right),
$$
$$
p_9=2 a_{99}-a_{91} \left( \varepsilon_9+\frac{1}{\varepsilon_1} \right)
-\frac{1}{\varepsilon_{59}} b_{59} X_5^{*}
\left( \frac{h_5}{h_9} \right)
-\frac{1}{\varepsilon_{10}} (b_{95} X_5^{*}+b_{10})
\left( \frac{h_{10}}{h_9} \right).
$$
Assuming
$$
\varepsilon_9=\frac{1}{\varepsilon_1},
$$
we select the value
$\varepsilon_1$
so that the equality is satisfied
$$
2 a_{11}-a_{19} \left( \varepsilon_1+\frac{1}{\varepsilon_9} \right)
=2 a_{99}-a_{91} \left( \varepsilon_9+\frac{1}{\varepsilon_1} \right).
$$
Then
$$
\varepsilon_1=\frac{(a_{11}-a_{99})
+\sqrt{(a_{11}-a_{99})^2+4a_{19} a_{91}}}{2a_{19}},
$$
and therefore,
$$
p_1=2\varepsilon
-\frac{1}{\varepsilon_2} \gamma_{21} M \left( \frac{h_2}{h_1} \right)
-\frac{1}{\varepsilon_{81}}
\gamma_{81} \frac{\rho_8}{\alpha_8} X_7^{*}
\left( \frac{h_8}{h_1} \right),
$$
$$
p_9=2\varepsilon
-\frac{1}{\varepsilon_{59}} b_{59} X_5^{*}
\left( \frac{h_5}{h_9} \right)
-\frac{1}{\varepsilon_{10}} (b_{95} X_5^{*}+b_{10})
\left( \frac{h_{10}}{h_9} \right),
$$
where
$\varepsilon$
is defined in~(2.5).
Thus, formula~(3.13) can be rewritten as
$$
U_{01}+W_{01}+W_{02}
=-2\delta \sum\limits_{j=1}^{10} h_j y_j^2(t)
-\sum\limits_{j=1}^{10} h_j (p_j-2\delta) y_j^2(t)
$$
$$
=-2\delta \sum\limits_{j=1}^{10} h_j y_j^2(t)
-\sum\limits_{j=1}^{10} h_j r_j y_j^2(t),
\eqno (3.14)
$$
where
$\delta$
is defined in~(2.6),
$$
r_1=2(\varepsilon-\delta)
-\frac{1}{\varepsilon_2} \gamma_{21} M \left( \frac{h_2}{h_1} \right)
-\frac{1}{\varepsilon_{81}}
\gamma_{81} \frac{\rho_8}{\alpha_8} X_7^{*}
\left( \frac{h_8}{h_1} \right),
$$
$$
r_2=2(\alpha_2-\delta)-\varepsilon_2 \gamma_{21} M
-\sum\limits_{k=3}^{7} \beta_k \left( \frac{h_k}{h_2} \right)
-\frac{1}{\varepsilon_3}
\left( b_{32}+b_3 X_5^{*} \right) X_3^{*}
\left( \frac{h_3}{h_2} \right)
$$
$$
-\frac{1}{\varepsilon_4} \left( b_{42}+b_4 X_6^{*} \right) X_4^{*}
\left( \frac{h_4}{h_2} \right)
-\frac{1}{\varepsilon_5} b_5 X_3^{*} X_5^{*}
\left( \frac{h_5}{h_2} \right)
-\frac{1}{\varepsilon_6} b_6 X_4^{*} X_6^{*}
\left( \frac{h_6}{h_2} \right),
$$
$$
r_3=2(\alpha_3-\delta)
-\frac{e^{\ae_3 \tau_3}}{\beta_3}
(b_{32} \rho_{32})^2 \Big( X_3^{*}+\theta_3 \Big)^2
-\varepsilon_3 \left( b_{32}+b_3 X_5^{*} \right) X_3^{*},
$$
$$
r_4=2(\alpha_4-\delta)
-\frac{e^{\ae_4 \tau_4}}{\beta_4}
(b_{42} \rho_{42})^2 \Big( X_4^{*}+\theta_4 \Big)^2
-\varepsilon_4 \left( b_{42}+b_4 X_6^{*} \right) X_4^{*},
$$
$$
r_5=2(\alpha_5-\delta)
-\frac{e^{\ae_5 \tau_5}}{\beta_5}
(b_5 \rho_5)^2
\Big( X_3^{*}+\theta_3 \Big)^2
\Big( X_5^{*}+\theta_5 \Big)^2
-\varepsilon_5 b_5 X_3^{*} X_5^{*}
-\varepsilon_{59} b_{59} X_5^{*},
$$
$$
r_6=2(\alpha_6-\delta)
-\frac{e^{\ae_6 \tau_6}}{\beta_6}
(b_6 \rho_6)^2
\Big( X_4^{*}+\theta_4 \Big)^2
\Big( X_6^{*}+\theta_6 \Big)^2
-\varepsilon_6 b_6 X_4^{*} X_6^{*},
$$
$$
r_7=2(\alpha_7-\delta)
-\frac{e^{\ae_7 \tau_7}}{\beta_7}
(b_7 \rho_7)^2
\Big( X_4^{*}+\theta_4 \Big)^2
\Big( X_6^{*}+\theta_6 \Big)^2
-\frac{1}{\varepsilon_{87}} \rho_8
\left( \frac{h_8}{h_7} \right),
$$
$$
r_8=2(\alpha_8-\delta)
-\varepsilon_{87} \rho_8
-\varepsilon_{81} \gamma_{81} \frac{\rho_8}{\alpha_8} X_7^{*},
$$
$$
r_9=2(\varepsilon-\delta)
-\frac{1}{\varepsilon_{59}} b_{59} X_5^{*}
\left( \frac{h_5}{h_9} \right)
-\frac{1}{\varepsilon_{10}} (b_{95} X_5^{*}+b_{10})
\left( \frac{h_{10}}{h_9} \right),
$$
$$
r_{10}=2(\alpha_{10}-\delta)
-\varepsilon_{10} (b_{95} X_5^{*}+b_{10}).
$$

We show that
$$
r_j \ge 0, \quad j=1,\dots,10.
$$
By virtue of definitions~(2.9)--(2.13) of the values
$\beta_k$,
$k=3,4,5,6,7$,
we have
$$
r_1=2(\varepsilon-\delta)
-\frac{1}{\varepsilon_2} \gamma_{21} M \left( \frac{h_2}{h_1} \right)
-\frac{1}{\varepsilon_{81}}
\gamma_{81} \frac{\rho_8}{\alpha_8} X_7^{*}
\left( \frac{h_8}{h_1} \right),
$$
$$
r_2=2(\alpha_2-\delta)-\varepsilon_2 \gamma_{21} M
-\frac{1}{\varepsilon_3}
\left(
b_{32} \rho_{32}
\Big( X_3^{*}+\theta_3 \Big) e^{\ae_3 \tau_3 /2}
+\left( b_{32}+b_3 X_5^{*} \right) X_3^{*}
\right)
\left( \frac{h_3}{h_2} \right)
$$
$$
-\frac{1}{\varepsilon_4}
\Big(
b_{42} \rho_{42}
\Big( X_4^{*}+\theta_4 \Big) e^{\ae_4 \tau_4 /2}
+\left( b_{42}+b_4 X_6^{*} \right) X_4^{*}
\Big)
\left( \frac{h_4}{h_2} \right)
$$
$$
-\frac{1}{\varepsilon_5}
\Big(
b_5 \rho_5
\Big( X_3^{*}+\theta_3 \Big)
\Big( X_5^{*}+\theta_5 \Big)
e^{\ae_5 \tau_5 /2}
+b_5 X_3^{*} X_5^{*}
\Big)
\left( \frac{h_5}{h_2} \right)
$$
$$
-\frac{1}{\varepsilon_6}
\Big(
b_6 \rho_6
\Big( X_4^{*}+\theta_4 \Big)
\Big( X_6^{*}+\theta_6 \Big)
e^{\ae_6 \tau_6 /2}
+b_6 X_4^{*} X_6^{*}
\Big)
\left( \frac{h_6}{h_2} \right)
$$
$$
-\frac{1}{\varepsilon_7}
b_7 \rho_7
\Big( X_4^{*}+\theta_4 \Big)
\Big( X_6^{*}+\theta_6 \Big)
e^{\ae_7 \tau_7 /2}
\left( \frac{h_7}{h_2} \right),
$$
$$
r_3=2(\alpha_3-\delta)
-\varepsilon_3
\left(
b_{32} \rho_{32} \Big( X_3^{*}+\theta_3 \Big) e^{\ae_3 \tau_3 /2}
+\left( b_{32}+b_3 X_5^{*} \right) X_3^{*}
\right),
$$
$$
r_4=2(\alpha_4-\delta)
-\varepsilon_4
\left(
b_{42} \rho_{42} \Big( X_4^{*}+\theta_4 \Big) e^{\ae_4 \tau_4 /2}
+\left( b_{42}+b_4 X_6^{*} \right) X_4^{*}
\right),
$$
$$
r_5=2(\alpha_5-\delta)
-\varepsilon_5
\left(
b_5 \rho_5
\Big( X_3^{*}+\theta_3 \Big)
\Big( X_5^{*}+\theta_5 \Big)
e^{\ae_5 \tau_5 /2}
+b_5 X_3^{*} X_5^{*}
\right)
-\varepsilon_{59} b_{59} X_5^{*},
$$
$$
r_6=2(\alpha_6-\delta)
-\varepsilon_6
\left(
b_6 \rho_6
\Big( X_4^{*}+\theta_4 \Big)
\Big( X_6^{*}+\theta_6 \Big)
e^{\ae_6 \tau_6 /2}
+b_6 X_4^{*} X_6^{*}
\right),
$$
$$
r_7=2(\alpha_7-\delta)
-\varepsilon_7
b_7 \rho_7
\Big( X_4^{*}+\theta_4 \Big)
\Big( X_6^{*}+\theta_6 \Big)
e^{\ae_7 \tau_7 /2}
-\frac{1}{\varepsilon_{87}} \rho_8
\left( \frac{h_8}{h_7} \right),
$$
$$
r_8=2(\alpha_8-\delta)
-\varepsilon_{87} \rho_8
-\varepsilon_{81} \gamma_{81} \frac{\rho_8}{\alpha_8} X_7^{*},
$$
$$
r_9=2(\varepsilon-\delta)
-\frac{1}{\varepsilon_{59}} b_{59} X_5^{*}
\left( \frac{h_5}{h_9} \right)
-\frac{1}{\varepsilon_{10}} (b_{95} X_5^{*}+b_{10})
\left( \frac{h_{10}}{h_9} \right),
$$
$$
r_{10}=2(\alpha_{10}-\delta)
-\varepsilon_{10} (b_{95} X_5^{*}+b_{10}).
$$
Considering definitions~(2.14)--(2.18) of the values
$\varepsilon_k$,
$k=3,4,5,6,7$,
and assuming
$$
\varepsilon_2=\frac{(\alpha_2-\delta)}{\gamma_{21} M},
$$
$$
\varepsilon_{59}=\frac{(\alpha_5-\delta)}{b_{59} X_5^{*}},
$$
$$
\varepsilon_{87}=\frac{(\alpha_8-\delta)}{\rho_8},
$$
$$
\varepsilon_{81}=\frac{(\alpha_8-\delta)}
{\gamma_{81} \frac{\rho_8}{\alpha_8} X_7^{*}},
$$
$$
\varepsilon_{10}=\frac{2(\alpha_{10}-\delta)}
{(b_{95} X_5^{*}+b_{10})},
$$
we obtain
$$
r_3=0,
\quad
r_4=0,
\quad
r_5=0,
\quad
r_6=0,
\quad
r_8=0,
\quad
r_{10}=0,
$$

$$
r_1=2(\varepsilon-\delta)
-\frac{1}{(\alpha_2-\delta)} (\gamma_{21} M)^2 \left( \frac{h_2}{h_1} \right)
-\frac{1}{(\alpha_8-\delta)}
\left( \gamma_{81} \frac{\rho_8}{\alpha_8} X_7^{*} \right)^2
\left( \frac{h_8}{h_1} \right),
$$
$$
r_2=(\alpha_2-\delta)
-\frac{1}{2(\alpha_3-\delta)}
\left(
b_{32} \rho_{32}
\Big( X_3^{*}+\theta_3 \Big) e^{\ae_3 \tau_3 /2}
+\left( b_{32}+b_3 X_5^{*} \right) X_3^{*}
\right)^2
\left( \frac{h_3}{h_2} \right)
$$
$$
-\frac{1}{2(\alpha_4-\delta)}
\Big(
b_{42} \rho_{42}
\Big( X_4^{*}+\theta_4 \Big) e^{\ae_4 \tau_4 /2}
+\left( b_{42}+b_4 X_6^{*} \right) X_4^{*}
\Big)^2
\left( \frac{h_4}{h_2} \right)
$$
$$
-\frac{1}{(\alpha_5-\delta)}
\Big(
b_5 \rho_5
\Big( X_3^{*}+\theta_3 \Big)
\Big( X_5^{*}+\theta_5 \Big)
e^{\ae_5 \tau_5 /2}
+b_5 X_3^{*} X_5^{*}
\Big)^2
\left( \frac{h_5}{h_2} \right)
$$
$$
-\frac{1}{2(\alpha_6-\delta)}
\Big(
b_6 \rho_6
\Big( X_4^{*}+\theta_4 \Big)
\Big( X_6^{*}+\theta_6 \Big)
e^{\ae_6 \tau_6 /2}
+b_6 X_4^{*} X_6^{*}
\Big)^2
\left( \frac{h_6}{h_2} \right)
$$
$$
-\frac{1}{(\alpha_7-\delta)}
(b_7 \rho_7)^2
\Big( X_4^{*}+\theta_4 \Big)^2
\Big( X_6^{*}+\theta_6 \Big)^2
e^{\ae_7 \tau_7}
\left( \frac{h_7}{h_2} \right),
$$
$$
r_7=(\alpha_7-\delta)
-\frac{1}{(\alpha_8-\delta)} \rho_8^2
\left( \frac{h_8}{h_7} \right),
$$
$$
r_9=2(\varepsilon-\delta)
-\frac{1}{(\alpha_5-\delta)} (b_{59} X_5^{*})^2
\left( \frac{h_5}{h_9} \right)
-\frac{1}{2(\alpha_{10}-\delta)} (b_{95} X_5^{*}+b_{10})^2
\left( \frac{h_{10}}{h_9} \right).
$$
Using definitions~(2.20)--(2.28) of the values
$h_j$,
$j=2,\dots,10$,
it follows that
$$
r_2=0,
\quad
r_7=0,
$$
$$
r_1=2(\varepsilon-\delta)
-\frac{1}{h_1}
\left[
\frac{5(\gamma_{21} M)^2}{(\alpha_2-\delta)^2}
+\left( \gamma_{81} \frac{\rho_8}{\alpha_8} X_7^{*} \right)^2
\frac{(\alpha_7-\delta)^2}
{(b_7 \rho_7 \rho_8)^2
\Big( X_4^{*}+\theta_4 \Big)^2
\Big( X_6^{*}+\theta_6 \Big)^2
e^{\ae_7 \tau_7}}
\right],
$$
$$
r_9=2(\varepsilon-\delta)
-\frac{1}{h_1} \frac{a_{91}}{a_{19}}
\left[
\frac{(b_{59} X_5^{*})^2}
{\Big(
b_5 \rho_5
\Big( X_3^{*}+\theta_3 \Big)
\Big( X_5^{*}+\theta_5 \Big)
e^{\ae_5 \tau_5 /2}
+b_5 X_3^{*} X_5^{*}
\Big)^2}
+1
\right].
$$
Finally, taking into account definition~(2.19) of the value
$h_1$
we get
$$
r_1 \ge 0,
\quad
r_9 \ge 0.
$$

So, from~(3.12) and~(3.14), it follows the estimate
$$
U_{01}+U_{02}+U_{\tau} \leq U_{01}+W_{01}+W_{02}+R_{\tau}
$$
$$
=-2\delta \sum\limits_{j=1}^{10} h_j y_j^2(t)
-\sum\limits_{j=1}^{10} h_j r_j y_j^2(t)+R_{\tau}
\leq -2\delta \sum\limits_{j=1}^{10} h_j y_j^2(t)+R_{\tau},
\eqno (3.15)
$$
where
$R_{\tau}$
is defined in~(3.8).



\subsection{The final estimate}

By virtue of estimates~(3.6) and~(3.15), from inequality~(3.1) we obtain
$$
\frac{d}{dt} V(t,y) \leq U_{01}+U_{02}+U_{03}+U_{\tau}
-\sum\limits_{k=3}^{7} \ae_k
\int\limits_{t-\tau_k}^{t} h_k \beta_k e^{-\ae_k (t-s)} y_2^2(s) ds
$$
$$
\leq -2\delta \sum\limits_{j=1}^{10} h_j y_j^2(t)
+q V^{3/2}(t,y)+R_{\tau}
-\sum\limits_{k=3}^{7} \ae_k
\int\limits_{t-\tau_k}^{t} h_k \beta_k e^{-\ae_k (t-s)} y_2^2(s) ds,
$$
where
$R_{\tau}$
is defined in~(3.8):
$$
R_{\tau}=h_3 \frac{e^{\ae_3 \tau_3}}{\beta_3}
(b_{32} \rho_{32})^2
\Big[
\Big( X_3^{*}+y_3(t-\tau_3) \Big)^2
-\Big( X_3^{*}+\theta_3 \Big)^2
\Big]
y_3^2(t)
$$
$$
+h_4 \frac{e^{\ae_4 \tau_4}}{\beta_4}
(b_{42} \rho_{42})^2
\Big[
\Big( X_4^{*}+y_4(t-\tau_4) \Big)^2
-\Big( X_4^{*}+\theta_4 \Big)^2
\Big]
y_4^2(t)
$$
$$
+h_5 \frac{e^{\ae_5 \tau_5}}{\beta_5}
(b_5 \rho_5)^2
\Big[
\Big( X_3^{*}+y_3(t-\tau_5) \Big)^2
\Big( X_5^{*}+y_5(t-\tau_5) \Big)^2
-\Big( X_3^{*}+\theta_3 \Big)^2
\Big( X_5^{*}+\theta_5 \Big)^2
\Big]
y_5^2(t)
$$
$$
+h_6 \frac{e^{\ae_6 \tau_6}}{\beta_6}
(b_6 \rho_6)^2
\Big[
\Big( X_4^{*}+y_4(t-\tau_6) \Big)^2
\Big( X_6^{*}+y_6(t-\tau_6) \Big)^2
-\Big( X_4^{*}+\theta_4 \Big)^2
\Big( X_6^{*}+\theta_6 \Big)^2
\Big]
y_6^2(t)
$$
$$
+h_7 \frac{e^{\ae_7 \tau_7}}{\beta_7}
(b_7 \rho_7)^2
\Big[
\Big( X_4^{*}+y_4(t-\tau_7) \Big)^2
\Big( X_6^{*}+y_6(t-\tau_7) \Big)^2
-\Big( X_4^{*}+\theta_4 \Big)^2
\Big( X_6^{*}+\theta_6 \Big)^2
\Big]
y_7^2(t).
$$
By virtue of designation~(2.29) of the value
$\omega$
and definition~(2.4) of the functional
$V(t,y)$,
from this, it follows the estimate
$$
\frac{d}{dt} V(t,y) \leq -2\omega V(t,y)+q V^{3/2}(t,y)+R_{\tau}.
\eqno (3.16)
$$

We introduce the notation
$$
\tau_{\min}=\min\{ \tau_3,\tau_4,\tau_5,\tau_6,\tau_7 \}>0.
$$
First, we assume that
$$
t \in (0,t') \cap (0,\tau_{\min}].
$$
Then
$t-\tau_k \leq 0$,
$k=3,4,5,6,7$,
and therefore,
$$
y_j(t-\tau_k)=\psi_j(t-\tau_k),
\quad j=3,4,5,6,
\quad k=3,4,5,6,7.
$$
Due to conditions~(2.31) and~(2.32)
$$
-X_j^{*} \leq \psi_j(t-\tau_k) \leq \theta_j,
\quad j=3,4,5,6,
\quad k=3,4,5,6,7,
$$
hence
$$
R_{\tau} \leq 0.
$$
Consequently, from inequality~(3.16) we have the estimate
$$
\frac{d}{dt} V(t,y) \leq -2\omega V(t,y)+q V^{3/2}(t,y).
\eqno (3.17)
$$
Since by virtue of condition~(2.31)
$\displaystyle
\sqrt{V(0,\psi)}<\frac{2\omega}{q}$,
then using the Gronwall inequality for
$t \in (0,t') \cap (0,\tau_{\min}]$
we obtain the estimate
$$
V(0,y) \leq \frac{V(0,\psi)}
{\displaystyle
\left( 1-\frac{q}{2\omega} \sqrt{V(0,\psi)} \right)^2}
\, e^{-2\omega t}.
$$
From this, if follows inequalities~(2.35):
$$
|y_j(t)| \leq \frac{1}{\sqrt{h_j}} \frac{\sqrt{V(0,\psi)}}
{\displaystyle
\left( 1-\frac{q}{2\omega} \sqrt{V(0,\psi)} \right)}
\, e^{-\omega t},
\quad
j=1,\dots,10.
$$
In particular, when
$j=10$
from condition~(2.33) the estimate follows
$$
|y_{10}(t)| \leq \frac{1}{\sqrt{h_{10}}} \frac{\sqrt{V(0,\psi)}}
{\displaystyle
\left( 1-\frac{q}{2\omega} \sqrt{V(0,\psi)} \right)}
\, e^{-\omega t}<e^{-\omega t}<1,
\quad t \in (0,t') \cap (0,\tau_{\min}].
$$
By virtue of the definition of the number
$t'$,
if
$t'<\infty$,
then
$y_{10}(t)<1$
for
$t \in (0,t')$
and
$y_{10}(t')=1$.
So,
$t'>\tau_{\min}$,
and for
$t \in (0,\tau_{\min}]$
inequalities~(2.34)--(2.35) are proved.

Let now
$$
t \in (0,t') \cap (\tau_{\min},2\tau_{\min}].
$$
Then
$t-\tau_k \leq \tau_{\min}$,
$k=3,4,5,6,7$.
If
$t-\tau_k \leq 0$,
then
$$
y_j(t-\tau_k)=\psi_j(t-\tau_k),
\quad j=3,4,5,6.
$$
If
$t-\tau_k \in (0,\tau_{\min}]$,
then for the functions
$y_j(t-\tau_k)$
the estimates of the form~(2.35) are valid:
$$
|y_j(t-\tau_k)| \leq \frac{1}{\sqrt{h_j}} \frac{\sqrt{V(0,\psi)}}
{\displaystyle
\left( 1-\frac{q}{2\omega} \sqrt{V(0,\psi)} \right)}
\, e^{-\omega (t-\tau_k)}
$$
$$
\leq \frac{1}{\sqrt{h_j}} \frac{\sqrt{V(0,\psi)}}
{\displaystyle
\left( 1-\frac{q}{2\omega} \sqrt{V(0,\psi)} \right)},
\quad j=3,4,5,6.
$$
In any case, due to conditions~(2.31)--(2.32) we have
$$
-X_j^{*} \leq y_j(t-\tau_k) \leq \theta_j,
\quad j=3,4,5,6,
\quad k=3,4,5,6,7,
$$
hence the inequality follows
$$
R_{\tau} \leq 0.
$$
Then from inequality~(3.16) we have estimate~(3.17).
Using the same reasoning as in the previous case,
we establish the validity of inequalities~(2.34)--(2.35) for
$t \in (\tau_{\min},2\tau_{\min}]$.

Using the method of mathematical induction,
a similar scheme proves the validity of inequalities~(2.34)--(2.35) for
$t \in (m\tau_{\min},(m+1)\tau_{\min}]$,
$m \in \mathbb{N}$.

So, inequalities~(2.34)--(2.35) are proved for all
$t>0$.

The theorem is completely proved.
\\

The author expresses the gratitude to Professor G.V.~Demidenko
for the attention to the research.

\end{document}